\documentclass[12pt,reqno]{amsart}
\usepackage[margin=0.95in]{geometry}
\usepackage{amsmath,amssymb,mathtools}
\usepackage{microtype}
\usepackage[hidelinks]{hyperref}
\hypersetup{
 pdftitle={Largest Element Orders Do Not Characterize Finite Simple Groups},
 pdfauthor={Jinbao Li, Wujie Shi, Yong Yang}
}

\newcommand{\PSp}{\operatorname{PSp}}
\newcommand{\PSL}{\operatorname{PSL}}
\newcommand{\Div}{\operatorname{Div}}
\newcommand{\vtwo}{v_2}
\newtheorem{theorem}{Theorem}[section]
\newtheorem{corollary}[theorem]{Corollary}
\newtheorem{lemma}[theorem]{Lemma}

\newcommand{\wt}{\lambda}

\title[Large element orders]{Largest Element Orders Do Not Characterize Finite Simple Groups}

\author[J.~B. Li]{Jinbao Li}
\address{School  of Mathematics and Physics, Suqian University, Jiangsu 223800, P. R. China}
\email{leejinbao25@hotmail.com}

\author[W.~J.  Shi]{Wujie Shi}
\address{School of Mathematical Sciences, Suzhou University, Jiangsu 215006, P. R. China;  and School of Mathematics and Artificial Intellegence,  Chongqing University of Arts and Sciences, Chongqing 402160, P. R. China}
\email{shiwujie@outlook.com}

\author[Y. Yang]{Yong Yang}
\address{Department of Mathematics, Texas State University, San Marcos, TX 78666, USA}
\email{yang@txstate.edu}


\begin{document}

\begin{abstract}
For every odd prime $p$ and every positive integer $k$, we construct two nonisomorphic nonabelian simple groups of Lie type in characteristic $p$ with the same order and the same $k$ largest distinct element orders. The groups are $\PSp_6(q)$ and $\Omega_7(q)$, where $q=p^{2^k}$. Thus the order together with any fixed number of largest element orders does not characterize all finite nonabelian simple groups.
\end{abstract}

\maketitle

\section{Introduction}

All groups considered in this paper are finite. For a group $G$, let
\[
 \omega(G)=\{|g|:g\in G\}
\]
denote its spectrum, and write its distinct members in decreasing order as
\[
 m_1(G)>m_2(G)>m_3(G)>\cdots.
\]
For a prime $p$, let
\[
 \omega_{p'}(G)=\{a\in\omega(G):p\nmid a\}
\]
denote its $p'$-spectrum.

A conjecture proposed by Shi in 1980s asserted that a  simple group $S$ is determined among finite groups by its order and spectrum, that is, if
\[
 |G|=|S|\qquad\text{and}\qquad \omega(G)=\omega(S),
\]
then $G\cong S$. This conjecture has been proved to
be true by Shi and Mazurov et al in a series of papers and is now a theorem; see \cite{Shi2023} and \cite{Shi2025} for a survey and references.

In this study, all elements of $\omega(S)$ of a simple group $S$ are involved. Hence, it is natural to seek weaker quantitative conditions under which a simple group $S$ can be reconstructed. A remarkable result of Kantor and Seress \cite{KS} shows that if $G$ and
$H$ are simple groups of Lie type of odd characteristic and
$m_i(G)=m_i(H)$ for $1\leq i\leq 3$, then $\operatorname{char}(G)=\operatorname{char}(H)$.
This means that the three largest element orders determine the
characteristic of simple groups of Lie type of odd characteristic.  Moreover, it was proved in \cite{He, ZS2013} that all  simple $K_3$-groups can be  characterized by using their orders and  largest and second largest element orders,  where a $K_n$-group is a group of order having  exactly $n$ distinct primes. It was  proved in \cite{GCC} that all sporadic simple groups can be determined by their orders and largest element orders. For $\PSL(2, p)$ with $p$ a prime, it was proved in \cite{LSY} that except $p=7$, every $\PSL(2, p)$ can be determined by its order and largest element order. All simple $K_4$-groups of type $\PSL(2, q)$ and simple $K_5$-groups of type $\PSL(3, p)$ with $(3, p-1)=1$ can be characterized in the same way by \cite{LYCS} and \cite{Yu} respectively.

Based on these observations, it was asked in \cite[Question~1]{Shi2018} whether one can replace the full spectrum by only $m_1$, by $\{m_1,m_2\}$, or by $\{m_1,m_2,m_3\}$ to characterize all nonabelian simple groups. Recall that by \cite[Theorem~1.2]{KS}, the three largest element orders determine the defining characteristic among simple groups of Lie type in odd characteristic.  However, this does not distinguish groups within the same characteristic. Indeed, $\PSp_{2n}(q)$ and $\Omega_{2n+1}(q)$, where $n\ge3$ and $q$ is odd and composite, have the same two largest element orders \cite[p.~804]{KS}.

We show that the agreement can extend to an arbitrarily long initial segment of the spectrum, even in any prescribed odd characteristic.

\begin{theorem}\label{thm:main}
Let $p$ be an odd prime, let $k$ be a positive integer, and set $q=p^{2^k}$. Then
\[
 S=\PSp_6(q)\qquad\text{and}\qquad G=\Omega_7(q)
\]
are nonisomorphic nonabelian simple groups such that
\[
 |S|=|G|
 \qquad\text{and}\qquad
 m_i(S)=m_i(G)\quad(1\le i\le k).
\]
\end{theorem}

Consequently, the order together with any fixed number of largest element orders does not characterize all nonabelian simple groups. This answers the question in \cite{Shi2018} negatively even when the comparison group is required to be simple.

The classification of nonisomorphic simple groups of the same order \cite[Theorem~5.1]{KLST} gives only the pairs $\PSp_{2n}(q)$ and $\Omega_{2n+1}(q)$, with $n\ge3$ and $q$ odd, and the pair $A_8$ and $\PSL_3(4)$. The latter pair has largest element orders $15$ and $7$, respectively \cite{Atlas}. Thus, when both groups are simple, counterexamples with equal largest element orders must come from the symplectic--orthogonal family. We use the rank-$3$ spectrum formulas in \cite[Corollaries~2 and~6]{Buturlakin} to show that the common semisimple orders can occupy an arbitrarily long initial segment.\\

The preceding discussion describes the situation when both comparison groups are simple. If the comparison group is allowed to be arbitrary, however, the symplectic--orthogonal pairs are not the only obstruction. The following example shows that the order and the three largest element orders can also fail to characterize an alternating simple group.

\begin{theorem}\label{thm:main+}
Let
\[
 S=A_{324},\qquad G=A_{323}\times C_2^2\times C_3^4.
\]
Then $S$ is nonabelian simple, $G\not\cong S$, and
\[
 |G|=|S|,\qquad m_i(G)=m_i(S)\quad(1\le i\le3).
\]
\end{theorem}

The smaller construction $S_{31}\times C_2^4$ has the same order and
largest element order as $A_{32}$; see \cite[Proposition~1.2]{MG}.
Theorem~\ref{thm:main+} shows that the analogous phenomenon can persist
through the three largest element orders. Its proof uses bounds on the
orders outside a common spectrum and does not require a computation of
either full spectrum. Further exact computations give a few additional
examples; these are recorded after the proof of Theorem~\ref{thm:main+}.

\section{Proof of Theorem \ref{thm:main}}

For a set $A$ of positive integers, let $\Div(A)$ be the set of positive divisors of its members. Write $[a_1,\ldots,a_s]$ for the least common multiple of positive integers $a_1,\ldots,a_s$, and let $\vtwo(a)$ denote the exponent of $2$ in $a$.

For  convenience, we first quote known results concerning the spectra of symplectic and orthogonal groups.

\begin{lemma} \label{lemma1}
Let $G=\PSp_{2n}(q)$, where $n\geq 2$ and $q$ is a power of odd prime $p$. Then $\omega(G)$ contains all divisors of the following numbers:

{\rm (1)} $\frac{q^n\pm 1}{2}$;

{\rm (2)} $[q^{n_1}+\epsilon_1 1, q^{n_2}+\epsilon_2 1, \ldots, q^{n_s}+\epsilon_s 1]$ for each $s\geq 2$, $\epsilon_i\in \{+, -\}$, $1\leq i\leq s$, and $n_1, n_2, \ldots, n_s>0$ such that $n_1+n_2+\cdots +n_s=n$;

{\rm (3)} $p^j\cdot [q^{n_1}+\epsilon_1 1, q^{n_2}+\epsilon_2 1, \ldots, q^{n_s}+\epsilon_s 1]$ for each $s\geq 1$, $\epsilon_i\in \{+, -\}$, $1\leq i\leq s$, and $j, n_1, n_2, \ldots, n_s>0$ such that $p^{j-1}+1+2n_1+2n_2+\cdots +2n_s=2n$;

{\rm (4)} $p^j$, if $p^{j-1}+1=2n$  for some $j>1$.
\end{lemma}

\begin{proof}
    This is \cite[Corollary 2]{Buturlakin}.
\end{proof}

\begin{lemma} \label{lemma2}
Let $G=\Omega_{2n+1}(q)$, where $n\geq 3$ and $q$ is a power of odd prime $p$. Then $\omega(G)$ contains all divisors of the following numbers:

{\rm (1)} $\frac{q^n\pm 1}{2}$;

{\rm (2)} $[q^{n_1}+\epsilon_1 1, q^{n_2}+\epsilon_2 1, \ldots, q^{n_s}+\epsilon_s 1]$ for each $s\geq 2$, $\epsilon_i\in \{+, -\}$, $1\leq i\leq s$, and $n_1, n_2, \ldots, n_s>0$ such that $n_1+n_2+\cdots +n_s=n$;

{\rm (3)} $p^j\cdot \frac{q^{n_1}\pm 1}{2}$ for all $j$ and $n_1$ such that $p^{j-1}+1+2n_1=2n$;

{\rm (4)} $p^j\cdot  [q^{n_1}+\epsilon_1 1, q^{n_2}+\epsilon_2 1, \ldots, q^{n_s}+\epsilon_s 1]$ for each $s\geq 2$, $\epsilon_i\in \{+, -\}$, $1\leq i\leq s$, and $j, n_1, n_2, \ldots, n_s>0$ such that $p^{j-1}+1+2n_1+2n_2+\cdots +2n_s=2n$;

{\rm (5)} $p^j$, if $p^{j-1}+1=2n$  for some $j>0$.
\end{lemma}

\begin{proof}
This is \cite[Corollary 6]{Buturlakin}. 
\end{proof}

\begin{proof}[Proof of Theorem~\ref{thm:main}]
Let $p,k,q,S,G$ be as in the theorem. Both groups are nonabelian simple, and
\begin{equation}\label{eq:order}
 |S|=|G|=\frac{q^9}{2}(q^2-1)(q^4-1)(q^6-1).
\end{equation}

Taking  $n=3$ in Lemma \ref{lemma1},  we obtain that $\omega(S)=\omega(PSp_6(q))$ consists of all divisors of the following numbers:

(1)  $\frac{q^3\pm 1}{2}$;

(2) $q\pm 1$, $[q\pm 1, q^2\pm 1]$;

(3) $9(q\pm 1)$ (if $p=3$), $p(q^2\pm 1)$, $p(q\pm 1)$; 

(4) $25$ (if $p=5$).

Letting $n=3$ in Lemma \ref{lemma2}, we have $\omega(G)=\omega(\Omega_7(q))$ consists of all divisors of the following numbers:

(1) $\frac{q^3\pm 1}{2}$;

(2) $q\pm 1$, $[q\pm 1, q^2\pm 1]$;

(3) $9\cdot \frac{q\pm 1}{2}$ (if $p=3$),  $p\cdot \frac{q^2\pm 1}{2}$;

(4)  $p(q\pm 1)$

(5) $25$ (if $p=5$).

Hence, the semisimple orders in either group are precisely the divisors of
\[
 \frac{q^3-1}{2},\quad \frac{q^3+1}{2},\quad q^2-1,\quad
 \frac{(q-1)(q^2+1)}2,\quad
 \frac{(q+1)(q^2+1)}2.
\]
In particular,
\begin{equation}\label{eq:pprime}
 \omega_{p'}(S)=\omega_{p'}(G).
\end{equation}
This equality holds for every odd prime power $q$ of characteristic $p$.

Since $q\ge p^2$ and $[q-1,q+1]=(q^2-1)/2$, we also have that 
\begin{equation}\label{eq:pmax}
 \begin{aligned}
 \max\{a\in\omega(S):p\mid a\}&=p(q^2+1),\\
 \max\{a\in\omega(G):p\mid a\}&=\frac{p(q^2+1)}2.
 \end{aligned}
\end{equation}
Indeed, when $p=3$ we have $q\ge9$ and $9(q+1)<3(q^2+1)$, while when $p=5$ we have $25<5(q^2+1)/2$. In particular, $p(q^2+1)\in\omega(S)\setminus\omega(G)$, so $S\not\cong G$.

Let
\[
 N=\frac{(q-1)(q^2+1)}2=[q-1,q^2+1].
\]
The semisimple formulas show that $N$ and all its divisors belong to the common $p'$-spectrum. Since $q=p^{2^k}$ is odd, $\vtwo(q^2+1)=1$, and the $2$-adic lifting-the-exponent formula gives
\[
 \vtwo(N)=\vtwo(q-1)
 =\vtwo(p-1)+k
 \ge k+1.
\]
Hence
\[
 N,\ \frac N2,\ \ldots,\ \frac{N}{2^{k-1}}
\]
are $k$ distinct common element orders.

Since $2^k\ge k+1$ and $p^k\ge2^k+1$, we have
\[
 q-1=p^{2^k}-1
 \ge p^{k+1}-1
 \ge p(2^k+1)-1
 >p\,2^k.
\]
Therefore
\[
 \frac{N}{2^{k-1}}
 =\frac{(q-1)(q^2+1)}{2^k}
 >p(q^2+1).
\]
By \eqref{eq:pmax}, all $k$ displayed common orders exceed every order divisible by $p$ in either group. Thus the $k$ largest element orders in both groups belong to their common $p'$-spectrum. Their equality follows from \eqref{eq:pprime}; together with \eqref{eq:order} and $S\not\cong G$, this proves the theorem.
\end{proof}

The smallest field in the construction already gives agreement of the four largest element orders.

\begin{corollary}\label{cor:q9}
Let $S=\PSp_6(9)$ and $G=\Omega_7(9)$. These are nonisomorphic nonabelian simple groups of the same order, and
\[
 (m_1,m_2,m_3,m_4)=(410,365,364,328)
\]
for both groups. Their fifth largest element orders are
\[
 m_5(S)=246\qquad\text{and}\qquad m_5(G)=205.
\]
\end{corollary}

\begin{proof}
Specializing Lemmas \ref{lemma1} and \ref{lemma2} to $q=9$ gives
\[
 \omega(S)=\Div\{410,365,364,328,246,240,90,72\}
\]
and
\[
 \omega(G)=\Div\{410,365,364,328,123,120,80,45,36\}.
\]
Every proper divisor of the first four displayed numbers is at most $205$, and all other generators are less than $328$. Thus the four largest orders are $410,365,364,328$ in both groups. The largest remaining order in $S$ is $246$, whereas that in $G$ is $205=410/2$.
\end{proof}

The equality of the semisimple spectra also gives the following consequence.

\begin{corollary}
For every odd prime $p$, the order together with the $p'$-spectrum does not characterize nonabelian simple groups, even among simple groups of Lie type in characteristic $p$.
\end{corollary}

\begin{proof}
Take $q=p^2$. Theorem~\ref{thm:main} gives nonisomorphic simple groups $\PSp_6(q)$ and $\Omega_7(q)$ of the same order, and \eqref{eq:pprime} gives equality of their $p'$-spectra.
\end{proof}

\section{Proof of Theorem \ref{thm:main+}}

For a positive integer $m$, define
\[
 \wt(m)=\sum_{p^a\parallel m}p^a,\qquad \wt(1)=0.
\]
The least degree of a permutation of order $m$ is $\wt(m)$.
The least degree of an even permutation of that order is
$\wt(m)$ when $m$ is odd, and $\wt(m)+2$ when $m$ is even.
Indeed, disjoint cycles of the indicated prime-power lengths attain the
first bound. If $m$ is even, a second even cycle is necessary for an
even permutation, and a disjoint transposition attains the additional
cost of two points; see \cite[Lemma 2.1]{MG} for example. 

All products over $p$ below are over primes. Define
\[
 \begin{aligned}
 B_1&=16\cdot27\cdot25\prod_{7\le p\le41}p,
       &\wt(B_1)&=296,\\
 B_2&=2\cdot9\cdot5\prod_{7\le p\le43}p,
       &\wt(B_2/2)&=285,\\
 B_3&=16\cdot25\prod_{7\le p\le43}p,
       &\wt(B_3)&=312.
 \end{aligned}
\]

\begin{lemma}\label{lem:bounds}
Let $m$ be a positive integer.
\begin{enumerate}
\item If $m$ has at most twelve distinct odd prime divisors, then
$m\le B_1\,2^{(\wt(m)-296)/12}$.
\item If $m$ is odd, then
$[m,6]\le B_2\,47^{(\wt(m)-285)/47}$.
\item If $3\nmid m$ and $m$ has at most twelve distinct odd prime
divisors, then $m\le B_3\,2^{(\wt(m)-312)/12}$.
\end{enumerate}
\end{lemma}

\begin{proof}
Consider first
\[
 f(x)=\frac{x^{12}}{2^x}\qquad(x>0).
\]
Among positive powers of $2$, $3$, and $5$, its largest values occur at
$16$, $27$, and $25$, respectively. This follows by comparing consecutive
powers, whose ratios are
\[
 \frac{f(p^{a+1})}{f(p^a)}
 =\frac{p^{12}}{2^{p^a(p-1)}}.
\]
For $p\ge7$, the inequality $2^{p-1}\ge p^2$ shows that this ratio is
less than one for every $a\ge1$. Thus the largest contribution from
such a prime is $f(p)$.

The function $f$ increases up to $12/\log2$ and then decreases.
The twelve largest odd-prime contributions are consequently attained at
\[
 27,25,7,11,13,17,19,23,29,31,37,41.
\]
The small-prime contributions exceed $f(41)$:
$f(25),f(27)>f(41)$ by monotonicity, while
$f(7)>f(41)$ follows from $7^{12}2^{34}>41^{12}$;
the contributions from $11$, $13$, and $17$ exceed $f(7)$.
All twelve listed values exceed one. If the prime $3$ is forbidden,
the twelve largest contributions are attained at
\[
 25,7,11,13,17,19,23,29,31,37,41,43.
\]
Here $f(7)>f(43)$ is equivalent to $56>43$ after taking twelfth
roots. Multiplying the contributions, and allowing the absent $2$-part
to contribute one, proves (1) and (3).

For (2), consider $h(x)=x^{47}/47^x$. If $p\ge5$, then the largest
value of $h$ on positive powers of $p$ occurs at $p$. For $p=5,7,11$,
the required comparisons are $5^{47}<47^{20}$, $7^{47}<47^{42}$,
and $11^{47}<47^{110}$; subsequent ratios decrease.
For $p\ge13$, the assertion follows from the
decrease of $h(x)$ for $x\ge13$. Moreover, for primes $p\ge5$,
$h(p)>1$ precisely when $p<47$, and $h(47)=1$.

The contribution of the prime $3$ to
$[m,3]^{47}/47^{\wt(m)}$ is
\[
 \begin{cases}
 3^{47},&3\nmid m,\\
 3^{47a}/47^{3^a},&3^a\parallel m.
 \end{cases}
\]
Its largest value occurs at $a=2$. In particular, this value exceeds
the contribution when $3\nmid m$, since $3^{47}>47^9$;
the comparison with $a=3$ is $3^{47}<47^{18}$.
The consecutive-power ratios settle all other exponents.
It follows that
\[
 \frac{[m,6]^{47}}{47^{\wt(m)}}
 \le \frac{B_2^{47}}{47^{285}},
\]
which proves (2).
\end{proof}

\begin{proof}[Proof of Theorem~\ref{thm:main+}]
Let $H=A_{323}$ and define
\[
 D=2^3\cdot3^2\cdot5\cdot43\prod_{7\le p\le37}p
   =3829100877561960,
 \qquad T=235D.
\]
The following three integers belong to $\omega(H)$:
\[
 \begin{array}{c|ccc}
 m&410D&246D&235D\\ \hline
 \wt(m)&321&319&319\\
 \wt(m)+2&323&321&321
 \end{array}
\]
They are distinct and at least $T$. We shall show that every order in
either $\omega(S)\setminus\omega(H)$ or
$\omega(G)\setminus\omega(H)$ is less than $T$.

Suppose first that $m\in\omega(A_{324})\setminus\omega(A_{323})$.
Then $\wt(m)=324$ if $m$ is odd, and $\wt(m)=322$ if $m$ is even.
In either case, the number of distinct odd prime divisors of $m$ is
even. It is at most twelve, since the first fourteen odd primes sum
to $326$. By Lemma \ref{lem:bounds}(1),
\[
 m\le B_1\,2^{(324-296)/12}
 =\frac{1230}{43}\,2^{7/3}D
 <\frac{9840}{43}D<T.
\]

The abelian factor $C_2^2\times C_3^4$ has exponent $6$.
Every element order in $G$ therefore divides $[m,6]$ for some
$m\in\omega(H)$. If $6\mid m$, the order already belongs to
$\omega(H)$.

If $m$ is odd, then $\wt(m)\le323$. By
Lemma \ref{lem:bounds}(2),
\[
 [m,6]
 \le B_2\,47^{38/47}
 =\frac{41}{4}\,47^{38/47}D<T.
\]
The last strict inequality is equivalent to the integer inequality
\[
 41^{47}47^{38}<940^{47}.
\]

If $m$ is even and $3\nmid m$, then $\wt(m)\le321$.
It has at most twelve distinct odd prime divisors, since the first
thirteen odd primes other than $3$ sum to $323$.
Lemma \ref{lem:bounds}(3) now gives
\[
 [m,6]=3m
 \le3B_3\,2^{9/12}
 =\frac{410}{3}\,2^{3/4}D<T.
\]
The final inequality is equivalent to
\[
 8\cdot82^4=361697408<395254161=141^4.
\]

Thus
\[
 \omega(S)\cap[T,\infty)
 =\omega(H)\cap[T,\infty)
 =\omega(G)\cap[T,\infty).
\]
There are at least three orders in these equal sets, so the three
largest distinct element orders of $S$ and $G$ coincide. Finally,
\[
 |G|=\frac{323!}{2}\,2^2 3^4
     =\frac{324!}{2}=|S|.
\]
The group $S$ is nonabelian simple, whereas $G$ has nontrivial center.
This proves the theorem.
\end{proof}


\medskip
An independent exact recurrence also gives
\[
 (m_1,m_2,m_3)=(410D,246D,235D)
\]
for both groups in Theorem~\ref{thm:main+}. The same computation gives
several further examples in which an alternating simple group and a
nonsimple group have the same order and the same three largest element
orders:
\[
\begin{array}{c|c}
 S & G\\ \hline
 A_{512} & S_{511}\times C_2^8\\
 A_{1296} & A_{1295}\times C_2^4\times C_3^4\\
 A_{7776} & S_{7775}\times C_2^4\times C_3^5\\
 A_{46656} & S_{46655}\times C_2^5\times C_3^6
\end{array}
\]
In each row, $|G|=|S|$ and $m_i(G)=m_i(S)$ for $1\le i\le3$.
We include these only to indicate that Theorem~\ref{thm:main+} is not an
isolated numerical occurrence; no general pattern is asserted here.

\section*{Acknowledgements}

The first author was  supported by the NSFC  (No. 12471017) and the Scientific Research Foundation for Advanced Talents of Suqian University (No. 2022XRC069). The third author was supported in part by Simons Foundation Award No.~918096. OpenAI's ChatGPT was used for mathematical exploration and editorial assistance. The authors checked the final arguments and are responsible for the contents of the paper.

\end{document}